\documentclass[12pt,a4paper]{article}

\usepackage[utf8]{inputenc}
\usepackage[english]{babel}
\usepackage{amsthm}
\usepackage{amsfonts}
\usepackage{amssymb}
\usepackage{amsmath}
\usepackage[affil-it]{authblk}
\usepackage{mathrsfs}
\usepackage{graphicx}
\usepackage{cite}

\newtheorem{Theorem}{Theorem}

\newtheorem{Lemma}{Lemma}
\newtheorem{Statement}{Statement}

\newenvironment{Def}{\par\bigskip\noindent{\bf Definition. }}{\par\bigskip}

\DeclareMathOperator{\re}{Re}
\DeclareMathOperator{\im}{Im}

\DeclareMathOperator{\Ai}{Ai}
\DeclareMathOperator{\Bi}{Bi}

\newcommand{\CC}{\mathbb{C}}
\newcommand{\RR}{\mathbb{R}}

\title{The Complex Airy Operator as Explicitly Solvable $\mathcal{PT}$--symmetrical Model
\thanks{The work is supported by Russian Science foundation, grant No 20-11-20261}}
\author[1]{Shkalikov A.\,A.\thanks{ashkaliko@yandex.ru}}
\author[2]{Tumanov S.\,N.\thanks{sergey.tumanov@yahoo.com}}
\affil[1]{Lomonosov Moscow State University}
\affil[2]{Moscow Center of Fundamental and Applied Mathematics at Lomonosov Moscow State University}

\begin{document}

\maketitle

\begin{abstract}
We study the Sturm--Liouville operator
$$
  T(\varepsilon)y=-\frac{1}{\varepsilon}y''+ p(x)y,
$$
with concrete $\mathcal{PT}$-- symmetric potential $p(x) = ix$ and Dirichlet boundary conditions on the segment $[-1,1]$. Here $\varepsilon \in (0, \infty)$ is a physical parameter. We explicitly describe a beautiful phenomenon of the eigenvalue behavior
when $\varepsilon$ changes from $0$  to $\infty$. All the critical values of $\varepsilon$  which determine the eigenvalue dynamics, are found in terms of the special Airy functions.
\end{abstract}


\section{Introduction}

The purpose of this work is to study the dynamics of the eigenvalues of the Sturm--Liouville operator
\begin{equation}\label{ix}
T(\varepsilon)y=-\frac{1}{\varepsilon}y''+ixy,\ \mathcal{D}_{T(\varepsilon)}=\bigl\{y\in W_2^2[-1,1]\,\bigl|\bigr.\, y(\pm1)=0\bigr\}.
\end{equation}
on the finite interval $[-1,1] $ with the parameter $\varepsilon >0$. For simplicity, the operator is considered
under the Dirichlet boundary conditions
\begin{equation}\label{Dir}
y(-1) = y(1) = 0.
\end{equation}
Actually, we also view in mind a more general problem:  to describe  the eigenvalue dynamics for the operator
$$
  T(\varepsilon)y=-\frac{1}{\varepsilon}y''+ p(x)y,
$$
with the so-called $\mathcal{PT}$--symmetrical potential $p(x) = \overline{p(-x)}$.
However, this problem is too difficult even for the case when $p(x)$ is a $\mathcal{PT}$--symmetrical polynomial. One can get in this case only partial information.  Namely, the spectrum of $T(\varepsilon)$ is real, provided that $\varepsilon \in (0, \varepsilon_0)$ for some $\varepsilon_0 = \varepsilon_0(p)$
(see \cite{TSh1}),
end for large $\varepsilon > R =R(\delta)$, where $\delta$ is arbitrary small, the spectrum lies in the $\delta$--neighborhood of the limit spectral graph, consisting
of some curves in the complex plane (the equations for these curves can be explicitly written down, see \cite{TSh2}, \cite{ES}).

In the general case we are not able to get an information on the behavior of the eigenvalues when $\varepsilon$ is changing in the middle interval $(\varepsilon_0, R(\delta))$. The exceptional case is  $p(x) =ix$. We will show that this case serves an
exactly solvable model. We will describe a beautiful phenomenon of the eigenvalue dynamics of the problem \eqref{ix}, \eqref{Dir} in terms of the special functions
$$
U_-(\xi) = -\sqrt 3 \Ai(\xi) + \Bi(\xi),\quad  U_+(\xi) = \sqrt 3 \Ai(\xi) + \Bi(\xi),
$$
where $ \Ai(\xi)$ and $ \Bi(\xi)$ are the standard Airy functions.

For each fixed $\varepsilon>0$ the operator $T(\varepsilon)$ has $\mathcal{PT}$--symmetry, that is, it
commutes with the composition $\mathcal{PT}$, where
$$
(\mathcal{T}y)(x)=\overline{y(x)},\quad (\mathcal{P}y)(x)=y(-x).
$$
The spectrum consists of real and pairwise conjugate complex eigenvalues.

A large number of works, especially in physics papers, are devoted to the study of $\mathcal{PT}$--symmetric operators.
Here we list some works \cite{BB1}--\cite{EG}, which stimulated our research.
More details on the literature can be found in the review by Dorey, Dunning and Tateo \cite{DDT2}.
Publications close to the topic of our study are mainly devoted to the proof of the reality of the spectrum of the one-dimensional
Schrodinger operator on the whole axis with some
special potentials, in particular, with cubic or polynomial potentials for certain values of the parameters. It was shown in \cite{GHH}, that the spectrum of the operator
\eqref{ix} on the whole line is empty. The detailed study of this operator and more general ones (when imaginary unit is changed by a complex number $c$)
 on the semi-axis was carried out in \cite{SS}. The properties of spectral projections for the complex valued potentials were studied in \cite {KS}  and \cite{MSV}.

In the papers \cite{Sh2}--\cite{ShT} the behavior of the spectrum with respect to the parameter $\varepsilon\to+\infty$ was studied in details. However, it was done only for large values of the parameter $\varepsilon$. Here we will get much more precise information even for large $\varepsilon$.
In particular, we will specify exactly the number  $\varepsilon_1$ such that for $0<\varepsilon\le\varepsilon_1$ the operator
$T(\varepsilon)$ has a purely real spectrum, as well as the critical numbers $\varepsilon_k$ which are the collision spectrum points with corresponding Jordan cells.

Numerical calculations show (they get their confirmation in the theorem \ref{MainTheorem}) that
starting from small values of the parameter $\varepsilon$ all eigenvalues of the operator family
$T(\varepsilon)$ move from infinity to the left. As $\varepsilon$ increases, at
$\varepsilon \approx 5.1$ the first eigenvalue crosses the knot point $1/\sqrt 3\approx 0.58$
and continues to move to the left. At $\varepsilon =\varepsilon_{1,\min} \approx 9.3$
the first eigenvalue reaches some point $\lambda_{1,\min}\approx 0.45<1/\sqrt 3$,
and as $\varepsilon >\varepsilon_{1,\min}$ increases further,
begins to move in the opposite direction, while all other eigenvalues continue to move to the left.
\begin{center}
\includegraphics[width=7cm,keepaspectratio]{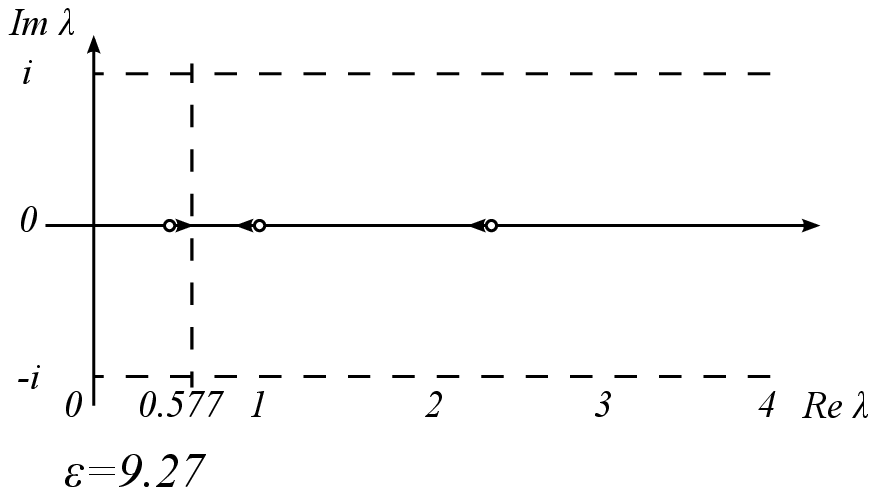}
\end{center}
Further, the first and second eigenvalues move towards each other,
approaching the knot point $1/\sqrt 3$, where they collide at $\varepsilon_1\approx12.3$.
\begin{center}
\includegraphics[width=7cm,keepaspectratio]{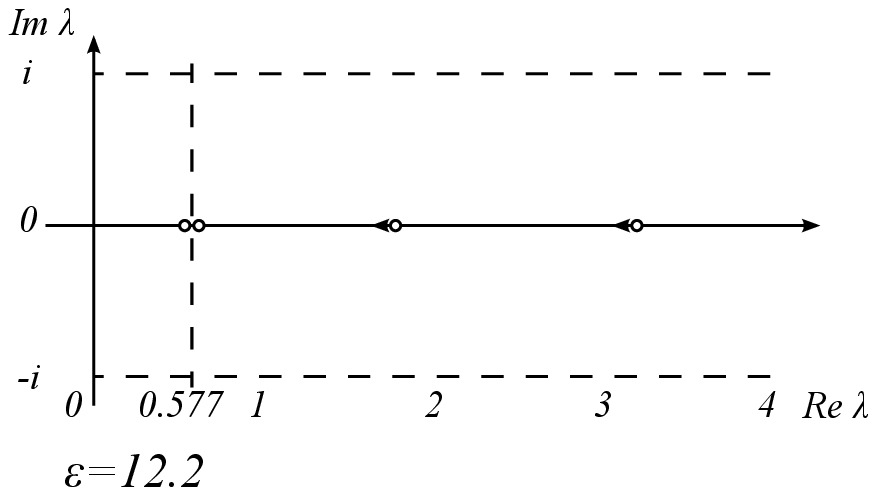}
\end{center}
After the collision at the knot point, the first and second eigenvalues move in the opposite direction to the complex plane perpendicularly to the real axis.
As $\varepsilon$ increases further, they approach the segments $[1/\sqrt 3, \pm i]$ and continue to move,
towards the points $\pm i$.
\begin{center}
\includegraphics[width=7cm,keepaspectratio]{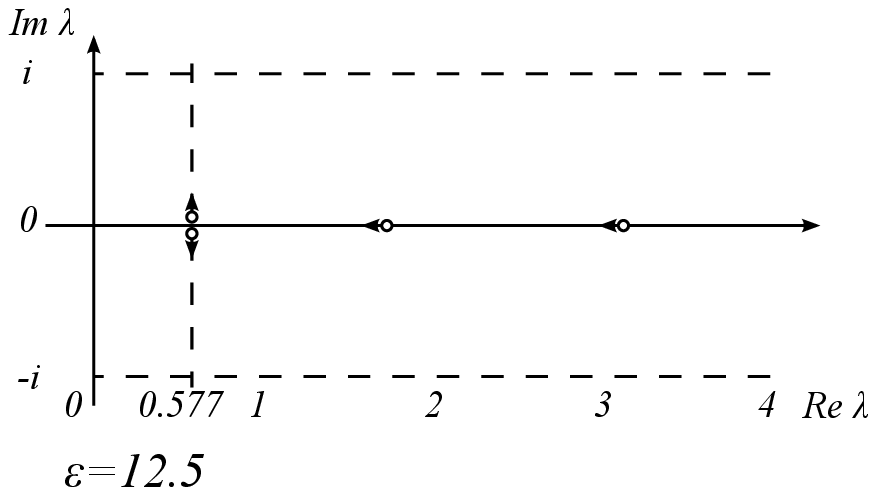}
\end{center}
Further, as $\varepsilon$ increase, the third eigenvalue crosses the knot point, moves to the left to the point $\lambda_{2,\min} < 1/\sqrt 3$,
and then moves back towards the fourth eigenvalue
untill the collision again at the knot point $1/\sqrt 3$ and the subsequent movement into the complex plane. Fifth and sixth eigenvalue (and subsequent $2k-1$-th and $2k$-th)
repeat the same dynamic. For large $\varepsilon$, the eigenvalues accumulate on the real ray $[1/\sqrt 3, +\infty)$,
and along the segments $[1/\sqrt 3, \pm i]$, moving to the points $\pm i$.
\begin{center}
\includegraphics[width=7cm,keepaspectratio]{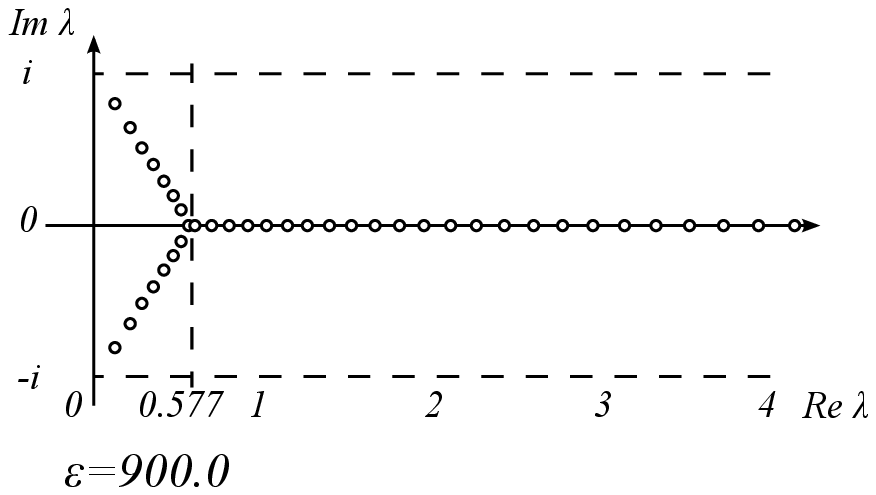}
\end{center}
Formulas for the asymptotic distribution of eigenvalues along the segments $[1/\sqrt 3, \pm i]$ and on the ray
$[1/\sqrt 3, +\infty)$ were found in \cite{Sh2}.

\section{Formulation of the main result and preliminaries}

As before we  work with the spectral problem:
\begin{equation}
\label{MainEq}
-\frac{1}{\varepsilon}y''+ixy=\lambda y,\quad y(-1)=y(1)=0,
\end{equation}
where $\varepsilon>0$ is the independent parameter, and $\lambda\in\CC$ is the spectral parameter.
\begin{Def}
The set $\mathcal{E}$ of pairs $(\varepsilon,\lambda)\in\mathbb{R_{+}}\times\mathbb{C}$, for which, for a given $\varepsilon>0$, the number $\lambda$ belongs to the spectrum
of \eqref{MainEq}, is called the {\it spectral locus} of the problem \eqref{MainEq}. The subset $\mathcal{E}_\mathbb{R}\subset\mathcal{E}$,
corresponding to real eigenvalues $\lambda$, is called the {\it real part} of the spectral locus of the problem \eqref{MainEq} or the {\it real
spectral locus}.
\end{Def}
The spectral locus of \eqref{MainEq} determines by the equation:
$D(\varepsilon,\lambda)=0$,
\begin{equation}
\label{Introduction_Eq_W}
D(\varepsilon,\lambda)=\left|
\begin{array}{cc}
y_1(-1) & y_1(1) \\
y_2(-1) & y_2(1)
\end{array}
\right|,
\end{equation}
where $y_1$ and $y_2$ are two arbitrary linearly independent solutions to the differential equation for fixed $\varepsilon>0$ and $\lambda\in\CC$.

The function $D$ is an analytic function of two arguments $\varepsilon\in\CC\setminus\{0\}$, $\lambda\in\CC$. The dynamics of eigenvalues as functions of the parameter
$\lambda_n=\lambda_n(\varepsilon)$ can be described with the implicit function theorem and the Weierstrass preparation theorem, to which we will turn later.
\begin{Def}
The point $(\varepsilon_0,\lambda_0)\in \mathcal{E}$ is called {\it critical} if $\partial D(\varepsilon_0,\lambda_0)/\partial\lambda=0$, i.e.
the implicit function $\lambda=\lambda(\varepsilon)$ has a singularity at $\varepsilon=\varepsilon_0$.
The values of $\varepsilon$ corresponding to the critical points of the spectral locus,
will also be called {\it critical}.
\end{Def}
Consider the classical Airy equation
$$
y''=\xi\cdot y, \quad y =y(\xi),
$$
and its two standard solutions, the functions $\Ai$ and $\Bi$
(see, for example,~\cite{Olver}; some of their properties are given in the Appendix to this article).
A remarkable role in what follows is played by special solutions to the Airy equation:
\begin{equation}
\label{eqUplusUminus}
U_-(\xi)=-\sqrt{3}\Ai(\xi)+\Bi(\xi),\quad
U_+(\xi)=\sqrt{3}\Ai(\xi)+\Bi(\xi).
\end{equation}

The following theorem is the main result of our paper.

\begin{Theorem}
\label{MainTheorem}
The zeros of $U_-$ and $U_+$ are located on the rays $\arg z=\pi/3 + 2\pi k/3$, $k=-1,0,1$ symmetrically about the origin.
Let $\{\alpha_k\}_{k=0}^\infty$ and $\{\beta_k\}_{k=1}^\infty$ be zeros of $U_-$ and $U_+$, respectively,
lying on the ray $l=\bigl\{z\,\bigl|\bigr.\arg z=\pi/3,\ |z|\ge0\bigr\}$, numbered in ascending order of absolute values. The zeros of both functions alternate:
$$
0=|\alpha_0|<|\beta_1|<|\alpha_1|<|\beta_2|<|\alpha_2|<|\beta_3|<\ldots.
$$
Let
$$
\delta_k=\Bigl(|\beta_k|\frac{\sqrt{3}}{2}\Bigr)^3,\quad
\varepsilon_k=\Bigl(|\alpha_k|\frac{\sqrt{3}}{2}\Bigr)^3,\quad k\ge1
$$
(obviously, $0<|\delta_1|<|\varepsilon_1|<|\delta_2|<|\varepsilon_2|<\ldots$). Then
\begin{itemize}
\item For $\varepsilon\in (0, \varepsilon_1)$ all eigenvalues of $T(\varepsilon)$ are real and simple. We numerate them in ascending order.
\item For $\varepsilon\in (0, \varepsilon_k)$, $k\ge1$, the eigenvalues of the pair $\lambda_{2k-1}(\varepsilon)$ and $\lambda_{2k}(\varepsilon)$ are real and simple.
At $\varepsilon = \delta_k$ the odd eigenvalue passes the knot point $\lambda_{2k-1}(\delta_k)=1/\sqrt{3}$. With firther growth of $\varepsilon$ it
reaches the minimum $\lambda_{2k-1,\, \min} < 1/\sqrt 3$, and turns back.
Further it collides at critical value $\varepsilon = \varepsilon_k$ with the even eigenvalue $\lambda_{2k-1}(\varepsilon_k)=\lambda_{2k}(\varepsilon_k)=1/\sqrt{3}$
at the knot point.
\item
After the collision, the eigenvalues move in the opposite direction to the complex plane perpendicularly to the real axis and
do not return to the real axis later. Outside the real axis, the eigenvalues do not collide.
\item
All critical points of the spectral locus coincide with the set
$$
\mathcal{M}=\Bigl\{
\left(\varepsilon_k,\frac{1}{\sqrt{3}}\right)
\Bigr\}_{k=1}^\infty.
$$
At critical values of $\varepsilon=\varepsilon_k$, the knot point $1/\sqrt{3}$ is a two-fold eigenvalue of
$T(\varepsilon_k)$, to which the Jordan cell corresponds.

\item At $\varepsilon=\delta_k$, the eigenfunctions for $\lambda_{2k-1}=1/\sqrt{3}$ are represented in an explicit form:
$$
y(z)=U_+\left(\delta_k^{1/3}\left(\frac{1}{\sqrt{3}}-iz\right)\right).
$$
\item At $\varepsilon=\varepsilon_k$, the eigenfunctions for $\lambda_{2k-1}=\lambda_{2k}=1/\sqrt{3}$ are represented in a closed form:
$$
y(z)=U_-\left(\varepsilon_k^{1/3}\left(\frac{1}{\sqrt{3}}-iz\right)\right).
$$

\item Let $\{z_k\}_{k=1}^\infty$ be the zeros of the function $\Bi$ in the first quarter of the complex plane numerated in increasing order of absolute values.
The following estimate is valid for points $\lambda_{2k-1,\min}$:
$$
\lambda_{2k-1,\min}<\cot\arg z_k<1/\sqrt{3}.
$$
\end{itemize}
\end{Theorem}
We replace the independent variable and spectral parameter in \eqref{MainEq}.

Consider $\eta=\varepsilon^{1/3}(\lambda-ix)$ as a new independent variable,
and $\xi=\varepsilon^{1/3}(\lambda+i)$ as a new spectral parameter (despite the fact that it will enter the boundary condition).

Define $w(\eta)=y(i(\eta/\varepsilon^{1/3}-\lambda))$.
As a result, the problem \eqref{MainEq} takes the form:
\begin{equation}
\label{Proves_Eq_Airy}
w''=\eta\cdot w,\quad w(\xi)=w(\xi-2i\varepsilon^{1/3})=0.
\end{equation}

The resulting equation, the Airy equation, will be considered in the entire complex plane. The $\xi$ parameter is an arbitrary complex number.

One-to-one correspondence
between spectra of \eqref{MainEq} and \eqref{Proves_Eq_Airy} is obvious. Thus when talking about the spectrum, we will not specify exactly what problem we are talking about.

The eigenvalues of the problem \eqref{Proves_Eq_Airy} are determined by the equation $F(\varepsilon,\xi)=0$:
\begin{equation}
\label{Introduction_Eq_F}
F(\varepsilon,\xi)=\left|
\begin{array}{cc}
w_1(\xi) & w_1(\xi-2i\varepsilon^{1/3}) \\
w_2(\xi) & w_2(\xi-2i\varepsilon^{1/3})
\end{array}
\right|,
\end{equation}
where $w_1$ and $w_2$ --- two independent solutions to the Airy equation.
\begin{Lemma}
\label{Lem_Xi_locus}
Given $\varepsilon>0$, the eigenvalues $\xi\in\CC$ of \eqref{Proves_Eq_Airy} are located
strictly in the I-st quarter of the complex plane (excluding real and imaginary axes).
\end{Lemma}
{\noindent\bf Proof.}
The spectrum of \eqref{MainEq} lies in the strip $\{\re\lambda>0,\ |\im\lambda|<1\}$ \cite{Sh3}. The Lemma follows
from the relation $\xi=\varepsilon^{1/3}(\lambda+i)$.\qquad$\Box$
\begin{Lemma}
\label{Proves_Lemma_partial}
The following statements concerning the partial derivatives of $D$ and $F$ are valid:
\begin{itemize}
\item
At the points of the spectrum, both partial derivatives $\partial D/\partial\varepsilon$ and $\partial D/\partial\lambda$
cannot vanish simultaneously.

\item
The partial derivative $\partial F/\partial\varepsilon$ is not vanishing at the points of the spectrum.

\item
At the points of the spectrum, the partial derivatives $\partial^k D/\partial\lambda^k$ and $\partial^k F/\partial\xi^k$, $k\ge1$
are vanishing or not vanishing simultaneously.
\end{itemize}
\end{Lemma}
{\noindent\bf Proof.} Let $(\varepsilon,\lambda)$ belong to the spectral locus of the problem \eqref{MainEq}, that is, $D(\varepsilon,\lambda)=0$.
We represent the defining $D$ solutions $y_j$, $j=1,2$ as follows:
$$
y_j(x)=w_j(\varepsilon^{1/3}(\lambda-ix)),
$$
where $w_j$ are two independent solutions to the Airy equation. To shorten the notation, we denote $\eta_\pm=\varepsilon^{1/3}(\lambda\mp i)$. Then
$$
D(\varepsilon,\lambda)=w_1(\eta_-)w_2(\eta_+)-w_1(\eta_+)w_2(\eta_-).
$$

Assuming that both partial derivatives are equal to zero, we obtain two equalities:
\begin{align*}
&w_1'(\eta_-)w_2(\eta_+)(\lambda+i)+w_1(\eta_-)w_2'(\eta_+)(\lambda-i)-\\
-& w_1'(\eta_+)w_2(\eta_-)(\lambda-i)-w_1(\eta_+)w_2'(\eta_-)(\lambda+i)=0.
\end{align*}
$$
w_1'(\eta_-)w_2(\eta_+)+w_1(\eta_-)w_2'(\eta_+)- w_1'(\eta_+)w_2(\eta_-)-w_1(\eta_+)w_2'(\eta_-)=0,
$$
whence:
\begin{equation}
\label{Proves_Eq_partial_w1w2}
w_1'(\eta_-)w_2(\eta_+)-w_1(\eta_+)w_2'(\eta_-)=0.
\end{equation}
Multiplying by $w_1(\eta_-)$, using further the fact that $D(\varepsilon,\lambda)=0$, we obtain:
$$
w_1(\eta_+)\Bigl\{
w_2(\eta_-)w_1'(\eta_-)-w_1(\eta_-)w_2'(\eta_-)
\Bigr\}=0.
$$

Since the curly brackets contain the Wronskian of linearly independent solutions to the Airy equation, then  $w_1(\eta_+)=0$.

Multiplying \eqref{Proves_Eq_partial_w1w2} by $w_2(\eta_-)$, carrying out similar calculations, we get $w_2(\eta_+)=0$,
which leads us to a contradiction: two independent solutions to a linear second-order differential equation cannot vanish at one point.
The first statement is proved.

We turn to $F$. Without limiting the generality let $w_1(\eta)=\Ai(\eta)$, $w_2(\eta)=\Bi(\eta)$. The equality $\partial F/\partial\varepsilon=0$ implies:
$$
\Ai(\xi)\Bi'(\xi-2\varepsilon^{1/3}i)-\Bi(\xi)\Ai'(\xi-2\varepsilon^{1/3}i)=0.
$$
Multiplying both sides by $\Ai(\xi-2\varepsilon^{1/3}i)$ and using the fact that at the points corresponding to the spectrum
$$
\Bi(\xi)\Ai(\xi-2\varepsilon^{1/3}i)=\Ai(\xi)\Bi(\xi-2\varepsilon^{1/3}i),
$$
we get:
$$
\Ai(\xi)\Bigl\{
\Ai(\xi-2\varepsilon^{1/3}i)\Bi'(\xi-2\varepsilon^{1/3}i)-\Bi(\xi-2\varepsilon^{1/3}i)\Ai'(\xi-2\varepsilon^{1/3}i)
\Bigr\}=0.
$$
The curly braces contain the Wronskian $\mathcal{W}(\Ai,\Bi)=1/\pi$ (statement \ref{Airy_Solutions_Statement_zeros}),
thus $\Ai(\xi)=0$, which leads us to a contradiction, since
the function $\Ai$ has no zeros in the I-st quarter of the complex plane (see \cite{Olver}) --- where $\xi$ lies according to the Lemma \ref{Lem_Xi_locus}.
The second statement of the Lemma is proved.

Note that the equality to zero of the partial derivatives $\partial^k D/\partial\lambda^k$ does not depend on the choice of independent solutions $y_j$.
Without loss of generality, let
$$
y_1(x)=\Ai(\varepsilon^{1/3}(\lambda-ix)),\quad y_2(x)=\Bi(\varepsilon^{1/3}(\lambda-ix)),
$$
and
$$
D(\varepsilon,\lambda)=\Ai(\varepsilon^{1/3}(\lambda+i))\Bi(\varepsilon^{1/3}(\lambda-i))-
\Bi(\varepsilon^{1/3}(\lambda+i))\Ai(\varepsilon^{1/3}(\lambda-i)),
$$
i.e.
$$
W(\varepsilon,\lambda)=F(\varepsilon,\xi(\lambda)),\mbox{ where }\xi(\lambda)=\varepsilon^{1/3}(\lambda+i),
$$
$$
\partial^k W(\varepsilon,\lambda)/\partial\lambda^k=\varepsilon^{k/3}\cdot\partial^k F(\varepsilon,\xi(\lambda))/\partial\xi^k,
$$
whence our statement follows. The Lemma is completely proved.\qquad$\Box$
\bigskip

The following Lemma is a direct corollary of the Weierstrass preparation theorem:
\begin{Lemma}
\label{StatmntIntrWeierstrass}
Let $W(\alpha,z)$ be an analytic function of two variables defined in the neighbor\-hood $|\alpha-\alpha_0|<\rho$, $|z-z_0|<r$ of the point
$(\alpha_0,z_0)\in\mathbb{C}^2$, with $\alpha_0\in\mathbb{R}$. In addition, $W(\alpha_0, z_0)=0$, $\partial^n W(\alpha_0,z_0)/\partial z^n=0$ as
$n=1,\ldots,k$ ($k\ge1$),
$\partial^{k+1} W(\alpha_0,z_0)/\partial z^{k+1}\ne0$, $\partial W(\alpha_0,z_0)/\partial\alpha\ne0$.

Then there are $0<\tilde r\le r$ and $0<\tilde\rho\le\rho$, that for any $\alpha$: $0<|\alpha-\alpha_0|<\tilde\rho$ the equation
$W(\alpha,z)=0$ has exactly $k+1$ roots $z=z_j(\alpha)$ in the neighbor\-hood $|z-z_0|<\tilde r$.

The $z_j$--images of intervals $(\alpha_0-\tilde\rho,\alpha_0)$ and $(\alpha_0,\alpha_0+\tilde\rho)$
are smooth pairwise disjoint Jordan arcs ($z_j$--arcs). When $\alpha$ changes from $\alpha_0-\tilde\rho$ to $\alpha_0$, the $z_j$--arcs have
semi-tangents at the point $z_0$ ---
collision rays emerging from $z_0$ at angles $2\pi/(k+1)$ with respect to each other; as $\alpha$ changes further from $\alpha_0$ to $\alpha_0+\tilde\rho$
the $z_j$--arcs also have semi-tangents at the point $z_0$ --- divergence rays directed along the bisectors of the angles formed by the collision rays.
\end{Lemma}
\section{Real spectral locus }
Due to the change of the spectral parameter $\xi=\varepsilon^{1/3}(\lambda+i)$ and the Lemma \ref{Lem_Xi_locus}, the real spectral
locus of the problem \eqref{MainEq} is described by
complex zeros of real solutions to the Airy equation, the ones lying in the I-st quarter of the complex $\xi$--plane.
\begin{Lemma}
\label{lm_xi_as_zeroz}
The pair $(\varepsilon,\lambda)\in\mathcal{E}_\mathbb{R}$ if and only if $\xi=\varepsilon^{1/3}(\lambda+i)$ is a zero
of some real solution to the Airy equation $V_a(\xi)$, $a\in\RR$.
\begin{equation}
\label{eqVaSolutions}
V_a(\eta)=a\Ai(\eta)+\Bi(\eta).
\end{equation}
\end{Lemma}
{\noindent\bf Proof.} Let $\lambda=\overline\lambda$ be the eigenvalue of \eqref{MainEq}, $y$ be the eigenfunction. Then for
$w(\eta)=y(i(\eta/\varepsilon^{1/3}-\lambda))$ we have $w(\xi)=w(\overline\xi)=0$.

By the Lemma \ref{Lem_Xi_locus}, the value $\xi$ lies in the I-st quarter, which means that the solution $w$ to the Airy equation can be represented as
$w(\eta)=C\cdot V_a(\eta)$ with some constant $C$, $a\in\CC$ ---
the solution $w$ certainly cannot be proportional to $\Ai$, since the latter has no zeros in the I-st quarter of the plane.

Since $V_a(\xi)=V_a(\overline\xi)=0$, then $V_{\overline a}(\xi)=\overline{V_a(\overline\xi)}=0$, hence $(a-\overline a)\Ai(\xi)=0$, i.e. $a\in\RR$.

Conversely, let for some  $a\in\RR$ and $\xi$ from the I-st quarter, $V_a(\xi)=0$. Then for $\varepsilon=(\im\xi)^3$ the function $V_a(\eta)$ will be the eigenfunction
of the problem \eqref{Proves_Eq_Airy},
and the corresponding $\lambda=\re\xi/\im\xi$ --- the eigenfunction of the main one \eqref{MainEq}.\qquad$\Box$
\begin{figure}[tb]
\begin{center}
\includegraphics[width=10cm,keepaspectratio]{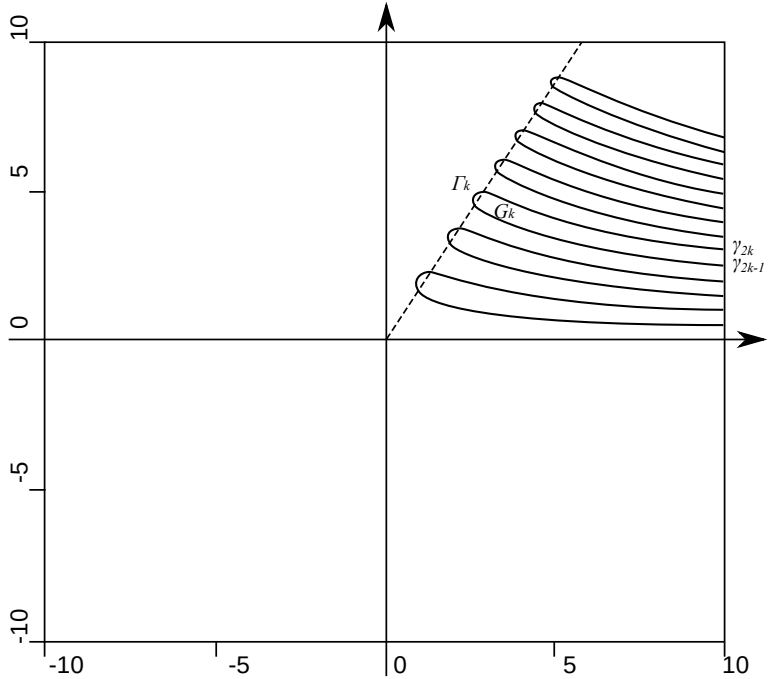}
\end{center}
\caption{Trajectories of zeros of real solutions to the Airy equation in the I-st quarter of the complex plane.}
\label{Proves_Pict_AiCurves}
\end{figure}
\begin{Lemma}
\label{Proves_Lemma_GammaN}
Given $a\in\mathbb{R}$ the equation $V_a(\xi)=0$ defines a countable number of implicit functions $\xi_n=\xi_n(a)$, $n\in\mathbb{Z}$ with values in the I-st quarter
of the complex plane, each continues analytically along the real axis.

The images $\Gamma_n$ of the real axis of each of these functions $\xi_n$ are pairwise disjoint Jordan analytic curves, unbounded from both
ends (see Fig.~\ref{Proves_Pict_AiCurves}).

The ends of $\Gamma_n$ asymptotically approache real positive semiaxis.
\end{Lemma}
{\noindent\bf Proof.} Fixing an arbitrary
real $a=a_0$, with Lemma~\ref{Airy_Solutions_Lemma_zeros} we find the infinite number of zeros of the function  $V_{a_0}$  in the I-st quarter. Denote them as
$\xi_n(a_0)$, $n\in\mathbb{N}$.

If $V_a(\xi)=0$, then $dV_a/d\xi\ne0$, since $V_a(\xi)$ is the nontrivial solution to the second order equation. We apply to $V(a,\xi)=V_a(\xi)$ the
implicit function theorem in neighborhoods of the points $(a_0,\xi_n(a_0))$ and find neighborhoods $\Upsilon_n(a_0)\subset\mathbb{C}$, in which
$\xi_n=\xi_n(a)$ will be single-valued analytic functions. Since for all $a\in\RR$ the derivative $dV_a/d\xi\ne0$, each of $\xi_n(a)$ can be uniquely continued
analytically along the entire real axis.

Thus, the curves $\Gamma_n=\{\xi_n(a)\,|\,a\in\RR\}$ are analytic.

Taking into account that $\Ai$ has no zeros in the I-st quarter, with the representation
$$
\frac{\Bi(\xi_n(a))}{\Ai(\xi_n(a))}=-a
$$
we obtain that $|\xi_n(a)|\to\infty$ as $a\to\pm\infty$, i.e. $\Gamma_n$ are unbounded from both
ends.

Self-intersections of $\Gamma_n$ would mean that two independent solutions $V_{a_1}(\xi)$ and $V_{a_2}(\xi)$ to the Airy equation have a common zero,
which is impossible.

The intersection of two different curves $\Gamma_n$ and $\Gamma_m$ is also impossible. Admitting the opposite --- intersection at some point
$\xi_n(a_1)=\xi_m(a_2)$,
--- we immediately obtain that $a_1=a_2$ (otherwise $V_{a_1}(\xi)$ and $V_{a_2}(\xi)$ are two independent solutions to the Airy equation with a common zero).
But then $\xi_n\equiv\xi_m$ by virtue of
the uniqueness of the analytic continuation and the uniqueness of the implicit function $\xi_n(a)$ satisfying the condition $V_{a_1}(\xi_n)=0$.

Expressing $\Bi(\xi)$ as a linear combination of $\Ai(\xi)$ and $\Ai(e^{-2\pi i/3}\xi)$ by the Statement \ref{Airy_Solutions_Statement_AiBi_Equiv},
we bring the equation $V_a(\xi)=0$ to the form:
$$
a+i=-2e^{-\pi i/6}\frac{\Ai(e^{-2\pi i/3}\xi)}{\Ai(\xi)}.
$$

Applying the asymptotics from the Statement \ref{Airy_Solutions_Statement_AiryAsymp} as $a\to\pm\infty$, $\xi\in\Gamma_n$ we obtain:
$$
a+i\sim-2\exp\Bigl(\frac{4}{3}\,\xi^{3/2}\Bigr).
$$
Separating the real and imaginary parts, we get $\tan(\im 4\xi^{3/2}/3)=O(1/a)$.

Let in the polar representation $\xi=re^{i\varphi}$, $r>0$, $0<\varphi<\pi/2$. Then for some
integer $m$
$$
r^{3/2}\sin\frac{3}{2}\varphi\sim\frac{3\pi}{4}m.
$$

As already noted $r\to\infty$, $0<\varphi<\pi/2$ as $a\to\pm\infty$; the value $m$ is constant due to analytical dependence $\xi(a)$.
We conclude that $\sin3\varphi/2=o(1)$, and $\varphi=o(1)$, thus
$$
r^{3/2}\frac{3}{2}\varphi\sim\frac{3\pi}{4}m,
$$
further, $\varphi=O(r^{-3/2})$, whence $\im\xi$=$r\sin\varphi=O(r^{-1/2})\to0$ as $a\to\pm\infty$. In other words, the ends of $\Gamma_n$
asymptotically approach the real positive semiaxis.\qquad$\Box$
\begin{Lemma}
\label{Proves_Lemma_Mu0Mu1_ZerosCritical}
Each of the curves $\Gamma_n$ intersects the ray $\{\arg\xi=\pi/3\}$ at least  at two points --- at one of the zeros $\alpha_k$ of the function $U_-$
and at one of the zeros $\beta_l$ of the function $U_+$.

The tangent to $\Gamma_n$ at the point $\alpha_k$ is parallel to the real axis,
and at $\beta_l$ is inclined at the angle $2\pi/3$ to the real axis.
\end{Lemma}
{\noindent\bf Proof.} Since each of the curves $\Gamma_n$ is parametrized as $\xi_n(a)$ by $a\in\RR$, we get with Lemma \ref{Airy_Solutions_Lemma_zeros}
$\xi_n(-\sqrt3)=\alpha_k$ --- some zero $U_-$, $\xi_n(\sqrt3)=\beta_l$ --- some zero $U_+$.

Differentiating the equality with respect to $a$:
$$
a\Ai(\xi_n(a))+\Bi(\xi_n(a))=0,
$$
taking into account the Wronskian $\mathcal{W}(\Ai,\Bi)=1/\pi$, we get:
\begin{equation}
\label{Proves_Eq_xiViaAi}
\xi_n'(a)=-\pi\Ai^2(\xi_n(a)).
\end{equation}

With formulas from the Statement~\ref{Airy_Solutions_Lemma_Mu0_Mu1} we express $\Ai$ in terms of $U_-$ and $U_+$:
$$
\Ai(z)\equiv\frac{U_+(z)-U_-(z)}{2\sqrt{3}}.
$$

Taking into account the symmetry of $U_-$ and $U_+$ following from the Statement \ref{Airy_Solutions_Lemma_Mu0_Mu1}, we obtain
$$
\Ai(\alpha_k)=\frac{U_+(\alpha_k)}{2\sqrt{3}}=\frac{U_+(-|\alpha_k|)}{2\sqrt{3}},\quad
\Ai(\beta_l)=-\frac{U_-(\beta_l)}{2\sqrt{3}}=e^{\pi i/3}\frac{U_-(-|\beta_l|)}{2\sqrt{3}},
$$
whence $\Ai^2(\alpha_k)>0$ and $\arg\Ai^2(\beta_l)=2\pi/3$. The Lemma follows with \eqref{Proves_Eq_xiViaAi}.\qquad$\Box$

\section{Proof of the main theorem  \ref{MainTheorem}}
Each $\Gamma_n$ splits the complex plane into two domains. One of them, denoted as $G_n$, has bounded
projection onto the imaginary axis. If follows from the Lemma \ref{Proves_Lemma_GammaN} as the ends of $\Gamma_n$ approach
the real positive semiaxis (see Fig. \ref{Proves_Pict_AiCurves}).

The formula $\xi=\varepsilon^{1/3}(\lambda+i)$ establishes a one-to-one correspondence between the spectral loci of the \eqref{MainEq} and \eqref{Proves_Eq_Airy} problems.
And due to Lemmas \ref{lm_xi_as_zeroz} and \ref{Proves_Lemma_GammaN} --- also between the real spectral locus $\mathcal{E}_\mathbb{R}$ of \eqref{MainEq} and
$\cup_n\Gamma_n$ in $\xi$--plane.

For real eigenvalues of \eqref{MainEq}, simple formulas are valid: $\varepsilon=(\im\xi)^3$,
$\lambda=\re\xi/\im\xi$.

The analytical dependence $\lambda_n=\lambda_n(\varepsilon)$ breaks at critical points $(\varepsilon,\lambda)$ of the
spectral locus --- where
$\partial D(\varepsilon,\lambda)/\partial\lambda=0$.
\begin{center}
***
\end{center}
{\bf Step 1.} {\it We prove that critical points $(\varepsilon,\lambda)$ of the real spectral locus \eqref{MainEq}
one-to-one correspond to the points of the curves $\Gamma_n$ with
tangents parallel to the real axis.}

If $(\varepsilon_0,\lambda_0)$ is a critical poin, denote $\xi_0=\varepsilon_0^{1/3}(\lambda_0+i)$.
Then $\xi_0\in\Gamma_n$ for some $\Gamma_n$ with
analytic parametrization $\xi_n(a)$ introduced by Lemma \ref{Proves_Lemma_GammaN}, $\xi_0=\xi_n(a_0)$, $a_0\in\RR$.

By Lemma \ref{Proves_Lemma_partial}, $\partial F(\varepsilon_0,\xi_0)/\partial\varepsilon\ne0$,
$\partial F(\varepsilon_0,\xi_0)/\partial\xi=0$. Applying the implicit function theorem, we represent locally
$\varepsilon=\varepsilon(\xi)$, $\varepsilon(\xi_0)=\varepsilon_0$,
$\varepsilon'(\xi_0)=0$. The function $\xi_n(a)$ is analytic in the neighborhood of the real axis. Differentiating the composition of analytic functions,
we get $d\varepsilon(\xi_n(a_0))/da=\varepsilon'(\xi_0)\xi_n'(a_0)=0$.

Along $\Gamma_n$ holds $\varepsilon(\xi_n(a))=(\im\xi_n(a))^3$, whence $\im\xi_n'(a_0)=0$ --- we use the reality of $a$ and the fact that $\im\xi_n(a)\ne0$ for all $a\in\RR$
(Lemma \ref{Lem_Xi_locus}).

Conversely, let $\im\xi_n'(a_0)=0$, $\xi_0=\xi_n(a_0)$. For $\varepsilon(\xi_n(a))=(\im\xi_n(a))^3$ we obtain that $d\varepsilon(\xi_n(a_0))/da=0$.

From \eqref{Proves_Eq_xiViaAi} we obtain that $\xi_n'(a)\ne0$ for all $a\in\RR$, as $\Ai(\xi)$ has no zeros in the I-st quarter of the complex plane.

Differentiating with respect to $a\in\RR$ the equality $F(\varepsilon(\xi_n(a)),\xi_n(a))\equiv0$, we get
$\partial F(\varepsilon_0,\xi_0)/\partial\xi=0$.

Lemma \ref{Proves_Lemma_partial} implies that $\partial D(\varepsilon_0,\lambda_0)/\partial\lambda=0$
for the corresponding eigenvalue  $\lambda_0$, i.e. $(\varepsilon_0,\lambda_0)$ is critical.
\begin{center}
***
\end{center}

For small $\varepsilon>0$ all eigenvalues of \eqref{MainEq} are real. We group them in pairs: the $k$-th pair consists
of $\lambda_{2k-1}$ and $\lambda_{2k}$, $k\ge1$.
We prove by induction on $k$ that
for $\varepsilon\in(0,\varepsilon_k)$ the eigenvalues with numbers $2k-1$ and higher are simple real, for
$\varepsilon=\delta_k<\varepsilon_k$ the eigenvalue $\lambda_{2k-1}(\delta_k)=1/\sqrt{3}$, for $\varepsilon=\varepsilon_k$ the eigenvalues of the $k$-th pair collide at
$\lambda_{2k-1}(\varepsilon_k)=\lambda_{2k}(\varepsilon_k)=1/\sqrt{3}$
and with further growth of $\varepsilon$ move in the opposite direction to the complex plane perpendicularly to the real axis, remain simple
and do not return back to the real axis later.
\begin{center}
***
\end{center}
{\bf Step 2 (The base case).} {\it We show that for $\varepsilon\in(0,\varepsilon_1)$ all eigenvalues are real and simple.

For $\varepsilon=\delta_1<\varepsilon_1$ the eigenvalue $\lambda_1(\delta_1)=1/\sqrt{3}$.

For $\varepsilon=\varepsilon_1$ the first pair collide at
$\lambda_1(\varepsilon_1)=\lambda_2(\varepsilon_1)=1/\sqrt{3}$, all other eigenvalues $\lambda_j(\varepsilon_1)$, $j>2$ remain real and simple.

With further growth of $\varepsilon$ the eigenvalues of the first pair move in the opposite direction to the complex plane perpendicularly to the real axis,
remain simple
and do not return back to the real axis later.}

It follows from the Lemma \ref{Proves_Lemma_Mu0Mu1_ZerosCritical}, there is at least one point on each $\Gamma_n$ corresponding to the critical point
of the spectral locus, with the tangent parallel to the real axis.

For small $\varepsilon$ all eigenvalues \eqref{MainEq} are simple and real (sufficient for $\varepsilon<3\pi^2/8$, \cite{TSh}), and remain so with
growth of $\varepsilon$ until it reaches the first critical
value, which we will denote as $\varepsilon=\hat\varepsilon$. We order all eigenvalues ascending for $\varepsilon<\hat\varepsilon$.
Let $\lambda_j$ be to the $j$-th eigenvalue, $j\in\mathbb{N}$.

For $0<\varepsilon<\hat\varepsilon$, each $\lambda_j$ corresponds to the analytic Jordan arc $\gamma_j$ in the first quarter of the $\xi$--plane.
The arcs are different and do not intersect.
Moreover, they are in one-to-one correspondence with the ends of the curves $\Gamma_n$.

Indeed, $\im\xi\to0$ as $\xi\to\infty$ along $\Gamma_n$ (Lemma \ref{Proves_Lemma_GammaN}), that is,
the ends of $\Gamma_n$ correspond to small $\varepsilon=(\im\xi)^3$.

Since all eigenvalues are real as $\varepsilon<\hat\varepsilon$,
and at critical points $\partial D/\partial\varepsilon\ne0$ (Lemma \ref{Proves_Lemma_partial}), taking into account the Lemma \ref{StatmntIntrWeierstrass}, we conclude
that the multiplicity of critical points for $\varepsilon=\hat\varepsilon$ is equal to two, i.e. for $\varepsilon=\hat\varepsilon$, the eigenvalues collide only in pairs,
after which they move in the opposite direction to the complex plane perpendicularly to the real axis --- a large multiplicity would
lead to a larger number of semi-tangents (collision rays), which
would be possible only in the presence of non-real eigenvalues.

Suppose $\lambda_j$ and $\lambda_{j+1}$ collide at $\varepsilon=\hat\varepsilon$. The union of  $\gamma_j$ and $\gamma_{j+1}$ completed by the point
$\hat\xi=\hat\varepsilon^{1/3}(\lambda_j+i)$, coincides with some $\Gamma_s$, since the curves $\Gamma_n$
are pairwise disjoint.

It follows from the construction that $\hat\xi$ is the only point on $\Gamma_s$ corresponding to the critical value of $\varepsilon$.
At the same time, $\Gamma_s$ intersects the ray $l=\bigl\{z\,\bigl|\bigr.\arg z=\pi/3\bigr\}$ at least at two points,
one of which ---  $\hat\beta$ the
zero of the function $U_{+}$, and the second --- $\hat\alpha$ --- the zero of the function $U_{-}$. And the last corresponds to the critical point of the locus.
In view of the uniqueness of such point on $\Gamma_s$, $\hat\alpha=\hat\xi$.

Since $\hat\varepsilon$ is the first critical point, it follows  $\hat\alpha=\alpha_1$ is the first non-real zero of $U_{-}$ on the $l$ ray.
If there were zeros on $l$ with smaller imaginary parts, due to the relation  $\varepsilon=(\im\xi)^3$, there would exist critical
values $\varepsilon<\hat\varepsilon$.

In other words, $\hat\varepsilon=\varepsilon_1$, $\hat\alpha=\alpha_1$. Since $\im\hat\beta<\im\alpha_1$, then $\hat\beta=\beta_1$ ---
the first zero of $U_{+}$ on $l$ corresponding to $\varepsilon=\delta_1<\varepsilon_1$.

For real eigenvalues $\lambda=\re\xi/\im\xi$. As $\beta_1\in l$, then
$\lambda_j(\delta_1)=1/\sqrt3$, i.e. at $\varepsilon=\delta_1$ one of the eigenvalues, which means that the minimal one $\lambda_1(\varepsilon)$, passes the point
$\tilde\lambda=1/\sqrt{3}$.
It follows that  $j=1$, the value $\varepsilon=\varepsilon_1$ corresponds to the collision of $\lambda_1$ and $\lambda_2$.
The rest of the eigenvalues remain simple and real.

As $\alpha_1\in l$, then $\lambda_{1}(\varepsilon_1)=\lambda_{2}(\varepsilon_1)=1/\sqrt3$.

We numerate the corresponding $\Gamma_s$ as $\Gamma_1$.

Since $\alpha_1$ corresponds to the critical $\varepsilon_1$, then $\partial F(\varepsilon_1,\alpha_1)/\partial\xi=0$,
$\partial F(\varepsilon_1,\alpha_1)/\partial\varepsilon\ne0$ (Lemma \ref{Proves_Lemma_partial}). As
the multiplicity of the first critical point is two,
$\partial^2 F(\varepsilon_1,\alpha_1)/\partial\xi^2\ne0$, with Lemma~\ref{StatmntIntrWeierstrass}
it follows that the divergence rays of the roots $\xi_1(\varepsilon)$, $\xi_2(\varepsilon)$ of the equation $F(\varepsilon,\xi)=0$ are directed along the normals to
$\Gamma_1$.

Therefore, for $\varepsilon>\varepsilon_1$ both roots $\xi_1(\varepsilon)$, $\xi_2(\varepsilon)$ fall into different areas,
bounded by the curve $\Gamma_1$. In particular, for some
$\delta>0$ for $\varepsilon_1<\varepsilon<\varepsilon_1+\delta$ there is exactly one root of the equation $F(\varepsilon,\xi)=0$ in the domain $G_1$, for definiteness,
$\xi_1(\varepsilon)$.

For all $\xi\in\Gamma_1$ we have $\varepsilon=(\im\xi)^3\le\varepsilon_1$, thus  for $\varepsilon>\varepsilon_1$, neither $\xi_1(\varepsilon)$ can cross $\Gamma_1$,
nor the others eigenvalues of \eqref{Proves_Eq_Airy}.

The $G_1$ domain does not contain points of $\Gamma_j$ for $j>1$, since all $\Gamma_j$ contain zeros of $U_{+}$ and $U_{-}$ on $l$,
and $G_1$ does not contain them, and $\Gamma_n$ themselves do not intersect in pairs, in particular, with $\Gamma_1$.

For finite $\varepsilon>\varepsilon_1$, $\xi_1(\varepsilon)$ cannot move to infinity, and other eigenvalues of \eqref{Proves_Eq_Airy} cannot come to $G_1$ from infinity.
This follows from the uniform in $\varepsilon>0$ boundedness of all non-real eigenvalues \cite{Sh3}, that is,
boundedness of the form $\xi=\varepsilon^{1/3}(\lambda+i)$ for finite  $\varepsilon$.

Thus, for $\varepsilon>\varepsilon_1$ the equation $F(\varepsilon,\xi)=0$ has a unique solution $\xi_1(\varepsilon)\in G_1$;
$\xi_1(\varepsilon)$ can be uniquely continued
analytically along the entire interval $(\varepsilon_1,+\infty)$, and the corresponding eigenvalue of \eqref{MainEq} in the $\lambda$--plane
does not return to the real axis and does not collide with other eigenvalues. The same statement is obviously true for the complex conjugate
eigenvalue.
\bigskip

{\bf Step 3 (The inductive step).} {\it Suppose that the statement holds for all indices $1,\cdots,k-1$.
We show its validity for the index $k$.}

For all $\varepsilon\le\varepsilon_{k-1}$, the eigenvalues with numbers $m\ge 2k-1$ remain simple and real. This is true until $\varepsilon$ reaches the next critical
$\varepsilon=\tilde\varepsilon>\varepsilon_{k-1}$. The eigenvalues with numbers less than
$2k-1$ do not return to the real axis, thus only pairs of real eigenvalues can collide, after which
they move in the opposite direction to the complex plane perpendicularly to the real axis.

Let $\lambda_j$ and $\lambda_{j+1}$ collide as $\varepsilon=\tilde\varepsilon$. The union of corresponding trajectories $\gamma_j$ and $\gamma_{j+1}$ in $\xi$--plane
completed by the point
$\tilde\xi=\tilde\varepsilon^{1/3}(\lambda_j+i)$, coincides with some $\Gamma_s$, intersecting with the ray $l$ at least
at $\tilde\beta$ the zero of $U_{+}$, and at $\tilde\alpha$ the zero of $U_{-}$. Like in the previous step, we conclude $\tilde\alpha=\tilde\xi$.

Due to the formula $\varepsilon=(\im\xi)^3$, we get $\tilde\xi=\tilde\alpha=\alpha_k$. Obviously $\tilde\alpha>\alpha_{k-1}$,
since for all $m<k$ $\alpha_m\in\Gamma_m$, on the other hand,
$\tilde\alpha\le\alpha_{k}$, as $\tilde\varepsilon$ is the next critical value following $\varepsilon_{k-1}$.

Likewise, $\tilde\beta=\beta_k$. As a result, $\tilde\varepsilon=\varepsilon_k$, the point $\delta_k=(\im\beta_k)^3$ corresponds to the
passage of the minimum of the remaining real
eigenvalues through the knot point $\tilde\lambda=1/\sqrt{3}$, i.e. $\lambda_{2k-1}(\delta_k)=1/\sqrt{3}$. Thus, $\lambda_j=\lambda_{2k-1}$,
$\lambda_{j+1}=\lambda_{2k}$. Further, $\lambda_{2k-1}(\varepsilon_k)=\lambda_{2k}(\varepsilon_k)=1/\sqrt{3}$.

We numerate the corresponding  $\Gamma_s$ as $\Gamma_k$.

Further reasoning showing that with increase of $\varepsilon$ the eigenvalues $\lambda_{2k-1}$ and $\lambda_{2k}$ do not return to the real axis and do not
collide with other eigenvalues outside the real axis, repeats the reasoning of the previous step.
\bigskip

{\bf Step 4.} {\it The estimate for $\lambda_{2k-1,\min}$.}

For $\varepsilon\in(\delta_k,\varepsilon_k)$, a certain part of the curve $\Gamma_k=\{\xi_k(a),\ a\in\RR\}$ corresponding to $a\in(-\infty,\sqrt{3}]$
match to the eigenvalue $\lambda_{2k-1}(\varepsilon)=\re\xi_k(a)/\im\xi_k(a)$,
$$
\lambda_{2k-1,\min}=\min_{\varepsilon\in(\delta_k,\varepsilon_k)}\lambda_{2k-1}(\varepsilon)\le\re\xi_k(0)/\im\xi_k(0)<1/\sqrt{3}.
$$

The Theorem follows with the note $\xi_k(0)=z_k$.\qquad$\Box$

\section{Appendix}
\setcounter{Statement}{0}
\setcounter{Lemma}{0}
\renewcommand{\theStatement}{A\arabic{Statement}}
\renewcommand{\theLemma}{A\arabic{Lemma}}
\nopagebreak

We recall some properties of the solutions to the classical Airy equation.
$$
y''=x\cdot y,
$$
It is known~\cite{Olver}, that the entire functions $\Ai$ and $\Bi$ form
the fundamental system of solutions (FSS) to this equation.
We present some Statements without proving. In more details
see \cite{Olver}.

\begin{Statement}
\label{Airy_Solutions_Statement_AiBi_Equiv}
For all complex $z$ the following identities hold:
$$
\Ai(z)+e^{2\pi i/3}\Ai(e^{2\pi i/3}z)+e^{-2\pi i/3}\Ai(e^{-2\pi i/3}z)\equiv 0,
$$
$$
\Bi(z)\equiv e^{\pi i/6}\Ai(e^{2\pi i/3}z)+e^{-\pi i/6}\Ai(e^{-2\pi i/3}z).
$$
\end{Statement}

\begin{Statement}
\label{Airy_Solutions_Statement_zeros}
The values of the Airy functions and their derivatives at zero are explicitly calculated:
$$
\Ai(0)=\frac{\Gamma(1/3)}{3^{1/6}2\pi}=\frac{1}{3^{2/3}\Gamma(2/3)}, \quad
\Ai'(0)=-\frac{3^{1/6}\Gamma(2/3)}{2\pi}=-\frac{1}{3^{1/3}\Gamma(1/3)},
$$
$$
\Bi(0)=\sqrt{3}\,\Ai(0)=\frac{1}{3^{1/6}\Gamma(2/3)}, \quad
\Bi'(0)=-\sqrt{3}\,\Ai'(0)=\frac{3^{1/6}}{\Gamma(1/3)}.
$$
The Wronskian
\begin{equation}
\label{Airy_Solutions_Eq_WronskianAiBi}
\mathcal{W}(\Ai,\Bi)=\Ai(0)\Bi'(0)-\Bi(0)\Ai'(0)=\frac{1}{\pi}.
\end{equation}
\end{Statement}

\begin{Statement}
\label{Airy_Solutions_Statement_AiryAsymp}
For any arbitrarily small  $\delta<\pi$ the asymptotic expansions are valid as $z\to\infty$ in the sector $|\arg z|\le\pi-\delta$:
$$
\Ai(z)\asymp\frac{e^{-\xi}}{2\pi^{1/2}z^{1/4}}\sum_{s=0}^\infty(-1)^s\frac{u_s}{\xi^s},\qquad
\Ai'(z)\asymp-\frac{z^{1/4}e^{-\xi}}{2\pi^{1/2}}\sum_{s=0}^\infty(-1)^s\frac{u_s}{\xi^s}.
$$

As $z\to\infty$ in the sector $|\arg z|\le\pi/3-\delta$:
$$
\Bi(z)\asymp\frac{e^{\xi}}{\pi^{1/2}z^{1/4}}\sum_{s=0}^\infty\frac{u_s}{\xi^s},\qquad
\Bi'(z)\asymp\frac{z^{1/4}e^{\xi}}{\pi^{1/2}}\sum_{s=0}^\infty\frac{u_s}{\xi^s}.
$$

In all formulas
$$
\xi=2z^{3/2}/3,\quad u_0=1,\quad u_s=\frac{2^s}{3^{3s}(2s)!}\frac{\Gamma\left(3s+\frac{1}{2}\right)}{\Gamma\left(\frac{1}{2}\right)}\mbox{ for } s\ge1,
$$
and fractional powers take principal values.
\end{Statement}

\begin{Statement}
\label{Airy_Solutions_via_varphi}
The function $\Ai$ has the following representation:
$$
\Ai(z)=\Ai(0)\varphi_1(z^3)+\Ai'(0)z\varphi_0(z^3),
$$
with entire functions $\varphi_0$, $\varphi_1$:
$$
\varphi_0(z)=1+\sum\limits_{k=1}^\infty\frac{\prod\limits_{s=0}^{k-1}\left(2+3 s\right)}{(1+3k)!}z^k,\quad
\varphi_1(z)=1+\sum\limits_{k=1}^\infty\frac{\prod\limits_{s=0}^{k-1}\left(1+3 s\right)}{(3k)!}z^k.
$$
\end{Statement}

For the solutions \eqref{eqUplusUminus} to the Airy equation the Statements \ref{Airy_Solutions_Statement_AiBi_Equiv}, \ref{Airy_Solutions_Statement_zeros} and
\ref{Airy_Solutions_via_varphi} lead to:
\begin{Statement}
\label{Airy_Solutions_Lemma_Mu0_Mu1}
$$
U_-(z)=\frac{2\cdot 3^{1/6}}{\Gamma(1/3)}z\varphi_0(z^3),\quad U_+(z)=\frac{2}{3^{1/6}\Gamma(2/3)}\varphi_1(z^3),
$$
$$
U_-(0)=0,\ U_-'(0)=\frac{2\cdot 3^{1/6}}{\Gamma(1/3)},\ U_+(0)=\frac{2}{3^{1/6}\Gamma(2/3)},\ U_+'(0)=0,
$$
$$
U_-(e^{\pm2\pi i/3}z)=e^{\pm2\pi i/3}U_-(z),\quad U_+(e^{\pm2\pi i/3}z)=U_+(z).
$$
\end{Statement}

Consider the family \eqref{eqVaSolutions} of real solutions to the Airy equation:
$$
V_a(z)=a\Ai(z)+\Bi(z),\ a\in\mathbb{R},
$$
describing, up to constant factors, all possible real solutions
with the exception of $\Ai$.

It is known~\cite{Olver} that the function $\Ai$ and its derivative
have zeros (infinite number, all are simple) only on the
negative real semiaxis.

We investigate the distribution of the zeros of $V_a(z)$ depending on the parameter $a\in\RR$.

\begin{Lemma}
\label{Airy_Solutions_Lemma_zeros}
The zeros of real solutions to the Airy equation and their derivatives are distributed as follows:

\begin{itemize}

\item As $a>\sqrt{3}$ the functions $V_a$ and $dV_a/dz$ have an infinite number of zeros on the negative semiaxis, an infinite number of
non-real zeros in the sector $\arg z\in(0,\pi/3)$ and their conjugates. Additionally, $dV_a/dz$ has exactly one zero on the positive semiaxis. There are no other zeros
of $V_a$ and $dV_a/dz$.

\item As $a=\sqrt{3}$, $V_a=U_+$, the functions $U_+$ and $dU_+/dz$  have an infinite number of zeros on the negative semiaxis, as well as zeros on the rays
$\arg z=\pm\pi/3$, obtained from the negative zeros by $\pm2\pi/3$ rotation of the complex plane. Additionally, $dU_+/dz$ has the 2-fold zero at $z=0$.
There are no other zeros of $U_+$  and $dU_+/dz$.

\item As $-\sqrt{3}<a<\sqrt{3}$ the functions $V_a$ and $dV_a/dz$ have an infinite number of zeros on the negative semiaxis, an infinite number of
non-real zeros in the sector $\arg z\in(\pi/3,\pi/2)$ and their conjugates. There are no other zeros of $V_a$ and $dV_a/dz$.

\item As $a=-\sqrt{3}$, $V_a=U_-$, the functions $U_-$ and $dU_-/dz$ have an infinite number of zeros on the negative semiaxis, as well as zeros on the rays
$\arg z=\pm\pi/3$, obtained from the negative zeros by $\pm2\pi/3$ rotation of the complex plane. Also, $U_-$ has zero at $z=0$.
There are no other zeros of $U_-$  and $dU_-/dz$.

\item As $a<-\sqrt{3}$ the functions $V_a$ and $dV_a/dz$ have an infinite number of zeros on the negative semiaxis, an infinite number of
non-real zeros in the sector $\arg z\in(0,\pi/3)$ and their conjugates. Additionally, $V_a$ has exactly one zero on the positive semiaxis. There are no other zeros
of $V_a$ and $dV_a/dz$.

\end{itemize}
All above mentioned zeros are simple with the exception of the 2-fold zero of $dU_+/dz$ at $z=0$.
Non real zeros approach the rays $\arg z=\pm\pi/3$ at infinity.
\end{Lemma}
{\noindent\bf Proof.} The zeros of real analytic functions $V_a$ and $dV_a/dz$ are either real or pairwise complex conjugate.

All real solutions to the Airy equation have an infinite number of
zeros on the negative semiaxis --- obviously follows from the existence of an infinite number of alternating zeros of $\Ai$ and $\Bi$ on the negative semiaxis \cite{Olver};

With the exception of $\Ai$, all real solutions to the Airy equation have infinitely many zeros in the right half-plane,
asymptotically approaching the rays $\arg z=\pm\pi/3$.
The statement follows from the asymptotic representations of the Statement~\ref{Airy_Solutions_Statement_AiryAsymp}.

All zeros of $V_a$ and $dV_a/dz$ are simple except for a 2-fold zero of $dU_+/dz$ at $z=0$
--- a consequence of the dimension of the space of solutions to a second-order equation. For the same reason, two solutions $V_{a_1}$ and $V_{a_2}$,
corresponding to different $a_1$ and $a_2$ cannot have common zeros. This also applies to the derivatives $dV_{a_1}/dz$ and $dV_{a_2}/dz$.

We investigate the distribution of non-real zeros by the Lommel method. For arbitrary complex constants $\alpha$ and $\beta$ and arbitrary
solution $w$ to the Airy equation one can directly verify following:
$$
\frac{d}{dx}\Bigl[\alpha w'(\alpha x) w(\beta x)-\beta w'(\beta x) w(\alpha x)
\Bigr]=
(\alpha^3-\beta^3)x w(\beta x) w(\alpha x),
$$
$$
\frac{d}{dx}\Bigl[\alpha^2 w'(\beta x) w(\alpha x)-\beta^2 w'(\alpha x) w(\beta x)
\Bigr]=
(\alpha^3-\beta^3) w'(\beta x) w'(\alpha x).
$$

Let $w$ be an arbitrary real solution to the Airy equation, $\alpha$ and $\beta$ two non-real conjugate zeros of $w$ or $w'$. We integrate these equalities
along the segment $[0,1]$, taking into account that $w(\alpha)=w(\beta)=0$ or $w'(\alpha)=w'(\beta)=0$. We get:
$$
-(\alpha-\beta)w(0) w'(0)=(\alpha^3-\beta^3)\int\limits_0^1 x |w(\alpha x)|^2\, dx,
$$
$$
-(\alpha^2-\beta^2)w(0) w'(0) = (\alpha^3-\beta^3) \int\limits_0^1 |w'(\alpha x)|^2\, dx.
$$

Assume that $\alpha=re^{\varphi i}$, $r>0$, $\arg\alpha=\varphi\in(0,\pi)$, $\beta=re^{-\varphi i}$. After minor transformations, we obtain:
\begin{align}
\label{Airy_Solutions_Eq_LemZer_1}
w(0)w'(0)\sin\varphi=-r^2\sin 3\varphi\int\limits_0^1 x |w(\alpha x)|^2\, dx,\\
\label{Airy_Solutions_Eq_LemZer_2}
w(0)w'(0)\sin 2\varphi=-r\sin 3\varphi\int\limits_0^1|w'(\alpha x)|^2\, dx.
\end{align}

Let $w(z)=V_a(z)$. Then
\begin{align*}
w(0)w'(0)
&=(a\Ai(0)+\Bi(0))(a\Ai'(0)+\Bi'(0))=\\
&=(a^2-3)\Ai(0)\Ai'(0)=-\frac{a^2-3}{2\pi\sqrt{3}}.
\end{align*}

Note that in \eqref{Airy_Solutions_Eq_LemZer_1} and \eqref{Airy_Solutions_Eq_LemZer_2} both $\sin\varphi$ and $\sin2\varphi$ are not equal to
zero: $\varphi\ne0$, $\varphi\ne\pi$, since
$\alpha$ is non-real zero,
and in case $\varphi=\pi/2$ the right-hand side in \eqref{Airy_Solutions_Eq_LemZer_2} is positive.

Further with \eqref{Airy_Solutions_Eq_LemZer_1} and \eqref{Airy_Solutions_Eq_LemZer_2} we get:
\begin{align*}
\frac{a^2-3}{2\pi\sqrt{3}}=r^2\frac{\sin 3\varphi}{\sin\varphi}\int\limits_0^1 x |w(\alpha x)|^2\, dx,\\
\frac{a^2-3}{2\pi\sqrt{3}}=r\frac{\sin 3\varphi}{\sin 2\varphi}\int\limits_0^1|w'(\alpha x)|^2\, dx.
\end{align*}

Analysis of the last two formulas shows:
\begin{itemize}
\item the value $|a|>\sqrt{3}$ if and only if $\varphi\in(0,\pi/3)$;
\item the value $|a|<\sqrt{3}$ if and only if $\varphi\in(\pi/3,\pi/2)$. In particular, this implies that non-real zeros of $\Bi$ ($a=0$)
lie in the region $|\arg z|\in(\pi/3,\pi/2)$;
\item the value $|a|=\sqrt{3}$ if and only if $\varphi=\pi/3$, i.e. up to a constant factor $U_-$ and $U_+$ are the only
solutions that may have zeros on the rays $\arg z=\pm\pi/3$.
\end{itemize}

None of the \eqref{eqVaSolutions} solutions have non-real zeros in the left half-plane or on the imaginary axis.

It remains show that $V_a$ (resp. $dV_a/dz$) has a single zero on the positive semiaxis if and only if $a<-\sqrt{3}$ (resp. $a>\sqrt{3}$).

The idea is the same for both cases; we limit ourselves by the case of the derivative. As we have already noted, $\Ai'$ has no zeros on the positive semiaxis.

Consider the real function
$$
v(z)=\frac{dV_a/dz}{\Ai'(z)}=a+\frac{\Bi'(z)}{\Ai'(z)}.
$$
According to \eqref{Airy_Solutions_Eq_WronskianAiBi} the derivative is negative on the positive semiaxis:
$$
v'(z)=\frac{\Bi''(z)\Ai'(z)-\Ai''(z)\Bi'(z)}{[\Ai'(z)]^2}=-\frac{z\mathcal{W}(\Ai,\Bi)}{[\Ai'(z)]^2}=-\frac{z}{\pi[\Ai'(z)]^2},
$$
and $v$ decreases for $z>0$. With the Statement~\ref{Airy_Solutions_Statement_zeros}, $v(0)=a-\sqrt{3}$.

It follows from the asymptotics of the Airy functions (the Statement~\ref{Airy_Solutions_Statement_AiryAsymp}) that $v(z)\to -\infty$ as $z\to+\infty$.
Thus $v$, and therefore $dV_a/dz$, has a unique positive zero if and only if $v(0)=a-\sqrt{3}>0$.
The lemma is completely proved. \qquad$\Box$

\end{document}